\documentclass[11pt]{article}

\usepackage{amsmath}
\usepackage{amssymb}
\usepackage{latexsym}

\def\be{\begin{equation}}
\def\ee{\end{equation}}
\def\bea{\begin{eqnarray}}
\def\eea{\end{eqnarray}}
\def\beas{\begin{eqnarray*}}
\def\eeas{\end{eqnarray*}}

\newtheorem{theorem}{Theorem}[section]
\newtheorem{definition}[theorem]{Definition}
\newtheorem{proposition}[theorem]{Proposition}
\newtheorem{corollary}[theorem]{Corollary}
\newtheorem{lemma}[theorem]{Lemma}
\newtheorem{remark}[theorem]{Remark}
\newtheorem{example}[theorem]{Example}
\newtheorem{examples}[theorem]{Examples}
\newtheorem{foo}[theorem]{Remarks}

\newenvironment{Remark}{\begin{remark}\rm}{\end{remark}}
\newenvironment{Remarks}{\begin{foo}\rm}{\end{foo}}
\newenvironment{proof}{\addvspace{\medskipamount}\par\noindent{\it Proof}.}
{\unskip\nobreak\hfill$\Box$\par\addvspace{\medskipamount}}

\newcommand{\ang}[1]{\left<#1\right>}  
\newcommand{\brak}[1]{\left(#1\right)}    
\newcommand{\crl}[1]{\left\{#1\right\}}   
\newcommand{\edg}[1]{\left[#1\right]}     

\newcommand{\E}[1]{{\rm E}\left[#1\right]}

\newcommand{\abs}[1]{\left|#1\right|}     

\def \ep{\hbox{ }\hfill$\Box$}

\begin{document}

\title{Second Order Backward Stochastic Differential Equations and Fully Non-Linear
Parabolic PDEs}

\author{Patrick Cheridito
        \thanks{Princeton University, Princeton USA,
        dito@princeton.edu.}
        \and
        H. Mete Soner
        \thanks{Ko\c{c} University, Istanbul, Turkey, msoner@ku.edu.tr. Member of the Turkish Academy of Sciences
                 and this work was partly supported by the Turkish Academy of Sciences.}
        \and
        Nizar Touzi
        \thanks{CREST, Paris, France, touzi@ensae.fr, and Tanaka
                Business School, Imperial College London, England, n.touzi@imperial.ac.uk.}
        \and Nicolas Victoir
        \thanks{Oxford University, Oxford, England, victoir@gmail.com.}
        }
\date{September 12, 2005}

\maketitle

\begin{abstract}
\noindent For a $d$-dimensional diffusion of the form $dX_t =
\mu(X_t) dt + \sigma(X_t) dW_t$, and continuous functions $f$ and
$g$, we study the existence and the uniqueness of adapted
processes $Y$, $Z$, $\Gamma$ and $A$ solving the second order
backward stochastic differential equation (2BSDE)
\beas
dY_t &=& f(t,X_t, Y_t, Z_t, \Gamma_t) \,dt + Z_t' \circ dX_t \, , \quad t \in [0,T) \, ,\\
dZ_t &=& A_t \,dt + \Gamma_t \,dX_t \, , \quad t \in [0,T) \, ,\\
Y_T &=& g(X_T) \, .
\eeas
If the associated PDE
\beas
&& - v_t(t,x) + f(t,x,v(t,x), Dv(t,x), D^2 v(t,x)) = 0 \, , \quad
(t,x) \in [0,T) \times \mathbb{R}^d \, ,\\
&& v(T,x) = g(x) \, ,
\eeas
admits a $C^3$-solution, it follows directly from It\^{o}'s lemma that
the processes
$$
v(t,X_t), \, Dv(t,X_t), \, D^2v(t,X_t), \, {\cal L} D v(t,X_t) \, ,
\quad t \in [0,T] \, ,
$$
solve the 2BSDE, where ${\cal L}$ is the Dynkin operator of $X$ without the drift term.

The main result of the paper shows that if the PDE has comparison
as in the theory of viscosity solutions and if $f$ is Lipschitz
in $Y$ and decreasing in $\Gamma$, the existence of a solution
$(Y,Z,\Gamma,A)$ to the 2BSDE implies that the associated PDE has
a unique continuous viscosity solution $v$, and $Y_t = v(t,X_t)$,
$t \in [0,T]$. In particular, the 2BSDE has at most one solution.
This provides a stochastic representation for solutions of fully
non-linear parabolic PDEs. As a consequence, the numerical
treatment of such PDE's can now be approached by Monte Carlo
methods.

\vspace{10pt}

 \vspace{10pt}

\end{abstract}


\setcounter{equation}{0}
\section{Introduction} \label{sectintroduction}

Since their introduction, backward stochastic differential
equations (BSDEs) have received considerable attention in the
probability literature. Interesting connections to partial
differential equations (PDEs) have been obtained and the theory
has found wide applications in areas like stochastic control,
theoretical economics and mathematical finance.

BSDEs were introduced by Bismut (1973) for the linear case and by
Pardoux and Peng (1990) for the general case. According to these
authors, a solution to a BSDE consists of a pair of adapted
processes $(Y,Z)$ taking values in $\mathbb{R}^n$ and
$\mathbb{R}^{d \times n}$, respectively, such that \be
\label{introbsde}
\begin{aligned}
dY_t &= f(t,Y_t,Z_t) dt + Z_t' dW_t \, , \quad t \in [0,T) \, ,\\
Y_T &= \xi \, ,
\end{aligned}
\ee
where $T$ is a finite time horizon, $(W_t)_{t \in [0,T]}$ a $d$-dimensional Brownian motion on
a filtered probability space $(\Omega, {\cal F}, ({\cal F}_t)_{t \in [0,T]}, P)$,
$f$ a progressively measurable function from
$\Omega \times [0,T] \times \mathbb{R}^n \times \mathbb{R}^{d \times n}$ to $\mathbb{R}^n$
and $\xi$ an $\mathbb{R}^n$-valued, ${\cal F}_T$-measurable random variable.

The key feature of BSDEs is the random terminal condition $\xi$
that the solution is required to satisfy. Due to the adaptedness
requirement on the processes $Y$ and $Z$, this condition poses
certain difficulties in the stochastic setting. But these
difficulties have been overcome, and now an impressive theory is
available; see for instance, Bismut (1973, 1978), Arkin and
Saksonov (1979), Pardoux and Peng (1990, 1992, 1994), Peng (1990,
1991, 1992a, 1992b, 1992c, 1993), Antonelli (1993), Ma, et al.
(1994, 1999, 2002), Douglas et al. (1996), Cvitani\'c and Ma (1996),
Cvitani\'c and Karatzas (1996), Chevance (1997), Cvitani\'c et al.
(1999), Pardoux and Tang (1999), Delarue (2002), Bouchard and
Touzi (2004), or the overview paper El Karoui et al. (1997).

If the randomness in the parameters $f$ and $\xi$ in
\eqref{introbsde} is coming from the state of a forward SDE, then
the BSDE is referred to as a forward-backward stochastic
differential equation (FBSDE) and its solution can be written as a
deterministic function of time and the state process. Under
suitable regularity assumptions, this function can be shown to be
the solution of a parabolic PDE. FBSDEs are called uncoupled if
the solution to the BSDE does not enter the dynamics of the
forward SDE and coupled if it does. The corresponding parabolic
PDE is semi-linear in case the FBSDE is uncoupled and quasi-linear
if the FBSDE is coupled; see Peng (1991, 1992b), Pardoux and Peng
(1992), Antonelli (1993), Ma et al. (1994), Pardoux and Tang
(1999), Ma and Yong (1999). These connections between FBSDEs and
PDEs have led to interesting stochastic representation results for
solutions of semi-linear and quasi-linear parabolic PDEs,
generalizing the Feynman--Kac representation of linear parabolic
PDEs and opening the way to Monte Carlo methods for the numerical
treatment of such PDEs, see for instance, Zhang (2001), Bally and
Pag\`es (2002), Ma et al (1994, 1999, 2002), Bouchard and Touzi
(2004), Delarue and Menozzi (2004). However, PDEs corresponding to
standard FBSDEs cannot be non-linear in the second order
derivatives because the second order terms only arise linearly
through It\^{o}'s formula from the quadratic variation of the
underlying state process.

In this paper we introduce FBSDEs with second order dependence in
the generator $f$. We call them second order backward stochastic
differential equations (2BSDEs) and show how they are related to
fully non-linear parabolic PDEs. This extends the range of
connections between stochastic equations and PDEs. In particular,
it opens the way for the development of Monte Carlo methods for
the numerical solution of fully non-linear parabolic PDEs. Our
approach is motivated by results in Cheridito et al. (2005a,
2005b)  which show how second order trading constraints lead to
non-linear parabolic Hamilton--Jacobi--Bellman equations for the
super-replication cost of European contingent claims.

The structure of the paper is as follows: In Section
\ref{sectnotation}, we explain the notation and introduce 2BSDEs
together with their associated PDEs. In Section
\ref{sectPDE2BSDE}, we show that the existence of a $C^3$ solution
to the associated PDE implies the existence of a solution to the
2BSDE. Our main result in Section \ref{sect2BSDEPDE} shows the
converse: If the PDE satisfies comparison as in the theory of
viscosity solutions and suitable Lipschitz and monotonicity
(parabolicity) conditions, then the existence of a solution to the
2BSDE implies that the PDE has a unique continuous viscosity
solution $v$. Moreover, the solution of the 2BSDE can then be
written in terms of $v$ and the underlying state process. This
implies that the solution of the 2BSDE is unique, and it provides
a stochastic representation result for fully non-linear parabolic
PDEs. In Section \ref{sectnumerical} we discuss Monte Carlo
schemes for the numerical solution of such PDEs. In Section 6 we
shortly discuss how the results of the paper can be adjusted to
the case of PDEs with boundary conditions.

\vspace{5mm}

\noindent {\bf Acknowledgements.}  Parts of this paper were
completed during a visit of Soner and Touzi to the Department of
Operations Research and Financial Engineering at Princeton
University. They would like to thank Princeton University, Erhan
\c{C}inlar and Patrick Cheridito for the hospitality. We also
thank Shige Peng and Arturo Kohatsu-Higa for stimulating and
helpful discussions.

\setcounter{equation}{0}
\section{Notation and definitions}
\label{sectnotation}

Let $d \ge 1$ be a natural number. We denote by
${\cal M}^d$ the set of all $d \times d$ matrices with
real components. $B'$ is the transpose of a matrix $B \in {\cal M}^d$ and
${\rm Tr}[B]$ its trace. By ${\cal M}^d_{\rm inv}$ we denote the set of all
invertible matrices in ${\cal M}^d$, by ${\cal S}^d$ all
symmetric matrices in ${\cal M}^d$, and by ${\cal S}^d_+$ all positive
semi-definite matrices in ${\cal M}^d$. For $B,C \in {\cal M}^d$,
we write $B \ge C$ if $B-C \in {\cal S}^d_+$. For $x \in \mathbb{R}^d$, we set
$$
|x| := \sqrt{x_1^2 + \dots x^2_d}
$$
and for $B \in {\cal M}^d$,
$$
|B| := \sup_{x \in \mathbb{R}^d \, , \,
|x| \le 1} B x \, .
$$
Equalities and inequalities between random variables are always
understood in the almost sure sense. We fix a finite time horizon
$T \in (0, \infty)$ and let $(W_t)_{t \in [0,T]}$ be a
$d$-dimensional Brownian motion on a complete probability space $(\Omega, {\cal F}, P)$.
For $t \in [0,T]$, we set
$W^t_s := W_s - W_t$, $s \in [t,T]$ and denote by
$\mathbb{F}^{t,T} = ({\cal F}^t_s)_{s \in [t,T]}$ the smallest filtration satisfying the
usual conditions and containing the filtration
generated by $(W^t_s)_{s \in [t,T]}$.

The coefficients $\mu : \mathbb{R}^d \to \mathbb{R}^d$ and $\sigma
: \mathbb{R}^d \to {\cal M}^d_{\rm inv}$ are continuous functions,
such that for all $N \ge 0$, there exists a constant $L_N$ with
 \bea \label{lip}
 |\mu(x) - \mu(y)| + |\sigma(x) - \sigma(y)|
 &\le&
 L_N |x-y|
 \quad \mbox{for all } x,y \in \mathbb{R}^d
 \mbox{ with } \abs{x}, \abs{y} \le N \, ,
 \eea
and constants $L \ge 0$ and $p_1 \in [0,1]$ such that \be
\label{growth} |\mu(x)| + |\sigma(x)| \le L(1 + |x|^{p_1}) \quad
\mbox{for all } x \in \mathbb{R}^d \, . \ee Then, for every
initial condition $(t,x) \in [0,T] \times \mathbb{R}^d$, the
forward SDE \be \label{sde}
\begin{aligned}
dX_s &= \mu(X_s) ds + \sigma(X_s) dW_s \, , \quad s \in [t,T] \, ,\\
X_t &= x \, ,
\end{aligned}
\ee has a unique strong solution $(X^{t,x}_s)_{s \in [t,T]}$. This
follows for instance, from Theorems 2.3, 2.4 and 3.1 in Chapter IV
of Ikeda and Watanabe (1989). Note that for existence and
uniqueness of a strong solution to the SDE \eqref{sde}, $p_1 =1$
in condition \eqref{growth} is enough. But for $p_1 \in [0,1)$, we
will get a better growth exponent $p$ in Proposition \ref{Vbound}
below. In any case, the process $(X^{t,x}_s)_{s \in [t,T]}$ is
adapted to the filtration $\mathbb{F}^{t,T}$, and by It\^{o}'s
formula we have for all $s \in [t,T]$ and $\varphi \in
C^{1,2}([t,T] \times \mathbb{R}^d)$,
 \bea \label{Ito}
 \varphi\left(s, X^{t,x}_s\right)
 &=&
 \varphi(t, x)
 + \int_t^s {\cal L} \varphi\left(r,X^{t,x}_r\right) dr
 + \int_t^s D \varphi\left(r, X^{t,x}_r\right)' dX^{t,x}_r \, ,
 \eea
where
 \beas
 {\cal L} \varphi (t,x) &=& \varphi_t(t,x)
 + \frac{1}{2} {\rm Tr} [D^2 \varphi(t,x) \sigma(x) \sigma(x)' ] \, ,
 \eeas
and $D \varphi$, $D^2 \varphi$ are the gradient and the matrix of second derivatives
of $\varphi$ with respect to the $x$ variables.

In the whole paper, $f : [0,T) \times \mathbb{R}^d \times
\mathbb{R} \times \mathbb{R}^d \times {\cal S}^d \to \mathbb{R}$
and $g : \mathbb{R}^d \to \mathbb{R}$ are continuous functions.

\begin{definition}
Let $(t,x) \in [0,T) \times \mathbb{R}^d$ and
$(Y_s,Z_s, \Gamma_s, A_s)_{s \in [t,T]}$ a quadruple of
$\mathbb{F}^{t,T}$-progressively measurable processes
taking values in $\mathbb{R}$, $\mathbb{R}^d$, ${\cal S}^d$ and $\mathbb{R}^d$,
respectively. Then we call $(Y,Z,\Gamma,A)$ a solution to the second order
backward stochastic differential equation {\rm (2BSDE)}
corresponding to $(X^{t,x},f,g)$ if
\bea
\label{2bsde1}
dY_s &=& f(s, X^{t,x}_s, Y_s ,Z_s , \Gamma_s) \,ds
+ Z_s' \circ dX^{t,x}_s \, , \quad s \in [t,T) \, ,\\
\label{2bsde2}
dZ_s &=& A_s \,ds + \Gamma_s \,dX^{t,x}_s \, , \quad s \in [t,T) \, ,\\
\label{2bsde3}
Y_T &=& g\left(X^{t,x}_T\right) \, ,
\eea
where $Z_s' \circ dX^{t,x}_s$ denotes Fisk--Stratonovich
integration, which is related to It\^{o} integration by
$$
Z_s' \circ dX^{t,x}_s =
Z_s' \,dX^{t,x}_s + \frac{1}{2} \,d\ang{Z, X^{t,x}}_s
=
Z_s' \,dX^{t,x}_s + \frac{1}{2} \,{\rm Tr}
[\Gamma_s \sigma(X^{t,x}_s) \sigma(X^{t,x}_s)' ] \,ds \, .
$$
\end{definition}
The equations \eqref{2bsde1}--\eqref{2bsde3} can be viewed as a whole
family of 2BSDEs indexed by
$(t,x) \in [0,T) \in \mathbb{R}^d$.
In the following sections, we will show relations between this family of
2BSDEs and the associated PDE
 \bea \label{pde}
 - v_t(t,x) + f\left(t,x,v(t,x),Dv(t,x),D^2v(t,x)\right)
 &=& 0 \quad \mbox{on } [0,T) \times \mathbb{R}^d \, ,
 \eea
with terminal condition
 \bea \label{terminal}
 v(T,x) &=& g(x) \, , \quad x \in \mathbb{R} \, .
 \eea
Since $Z$ is a semi-martingale, the use of the Fisk--Stratonovich integral in
\eqref{2bsde1} means no loss of generality, but it simplifies the
notation in the PDE \eqref{pde}.
Alternatively, \eqref{2bsde1} could be written in terms of
the It\^{o} integral as
 \bea \label{2bsdetilde}
 dY_s
 &=&
 \tilde{f}\left(s, X^{t,x}_s, Y_s ,Z_s , \Gamma_s\right) \,ds
 + Z_s' \,dX^{t,x}_s
 \eea
for
 \beas
 \tilde{f}(t,x,y,z,\gamma)
 &=&
 f(t,x,y,z,\gamma)
 + \frac{1}{2}\,{\rm Tr}\left[\gamma \sigma(x) \sigma(x)'\right]
 \, .
 \eeas
In terms of $\tilde{f}$, the PDE \eqref{pde} reads as follows:
 \beas
 - v_t(t,x) + \tilde{f}(t,x,v(t,x),Dv(t,x),D^2v(t,x))
 - \frac{1}{2} {\rm Tr} [D^2 v(t,x) \sigma(x) \sigma(x)']
 &=&
 0 \,.
 \eeas
Note that the form of the PDE \eqref{pde} does not depend on the
functions $\mu$ and $\sigma$ determining the dynamics in
\eqref{sde}. So, we could restrict our attention to the case where
$\mu \equiv 0$ and $\sigma \equiv I_d$, the
$d \times d$ identity matrix. But the
freedom to choose $\mu$ and $\sigma$ from a more general class
provides additional flexibility in the design of the Monte Carlo
schemes discussed in Section \ref{sectnumerical} below.

\setcounter{equation}{0}
\section{From a solution of the PDE to a solution of the 2BSDE}
\label{sectPDE2BSDE}

Assume $v : [0,T] \times \mathbb{R}^d \to \mathbb{R}$ is a
continuous function such that
$$
v_t, Dv, D^2v, {\cal L} D v \mbox{ exist and are continuous on }
[0,T) \times \mathbb{R}^d \, ,
$$
and $v$ solves the PDE \eqref{pde} with terminal condition \eqref{terminal}.
Then it follows directly from It\^{o}'s formula \eqref{Ito} that
for each pair $(t,x) \in [0,T) \times \mathbb{R}^d$, the
processes
 \beas
 Y_s &=& v\left(s,X^{t,x}_s\right) \, , \quad s \in [t,T] \, ,\\
 Z_s &=& Dv\left(s, X^{t,x}_s\right) \, , \quad s \in [t,T] \, ,\\
 \Gamma_s &=& D^2v\left(s, X^{t,x}_s\right) \, , \quad s \in [t,T] \, ,\\
 A_s &=& {\cal L} Dv \left(s, X^{t,x}_s\right) \, , \quad s \in [t,T] \, ,
 \eeas
solve the 2BSDE corresponding to $(X^{t,x}, f,g)$.

\setcounter{equation}{0}
\section{From a solution of the 2BSDE to a solution of the PDE}
\label{sect2BSDEPDE}

In all of Section \ref{sect2BSDEPDE} we assume that
$$
f : [0,T) \times \mathbb{R}^d \times \mathbb{R} \times
\mathbb{R}^d \times {\cal S}^d \to \mathbb{R} \quad \mbox{and}
\quad g : \mathbb{R}^d \to \mathbb{R} \,
$$
are continuous functions that satisfy the following Lipschitz and growth assumptions:

{\bf (A1)} \quad For every $N \ge 1$ there exists a constant $F_N$
such that
 \beas
 \abs{f(t,x,y,z,\gamma) - f(t,x,\tilde{y}, z, \gamma)}
 &\le&
 F_N \abs{y - \tilde{y}}
 \eeas
for all $(t,x, y, \tilde{y} ,z,\gamma) \in [0,T] \times
\mathbb{R}^d \times \mathbb{R}^2 \times \mathbb{R}^d \times {\cal
S}^d$ with $\max \crl{\abs{x} \, , \, \abs{y} \, , \,
\abs{\tilde{y}} \, , \, \abs{z} \, , \, \abs{\gamma}} \le N$.

{\bf (A2)} \quad There exist constants $F$ and $p_2 \ge 0$ such
that
$$
|f(t,x,y,z,\gamma)| \le F (1 + |x|^{p_2} + |y| + |z|^{p_2} + |\gamma|^{p_2})
$$
for all $(t,x,y,z,\gamma) \in [0,T] \times \mathbb{R}^d \times
\mathbb{R} \times \mathbb{R}^d \times {\cal S}^d$.

{\bf (A3)} \quad There exist constants $G$ and $p_3 \ge 0$ such
that
 \beas
 |g(x)|
 &\le&
 G (1 + |x|^{p_3})
 \quad \mbox{for all } x \in \mathbb{R}^d
 \, .
 \eeas

 \subsection{Admissible strategies}

We fix constants $p_4, p_5 \ge 0$ and denote for all
$(t,x) \in [0,T] \times \mathbb{R}^d$ and $m \ge 0$ by ${\cal A}_m^{t,x}$
the class of all processes of the form
$$
Z_s = z + \int_t^s A_r dr + \int_t^s \Gamma_r dX^{t,x}_r \, ,
\quad s \in [t,T] \, ,
$$
where $z \in \mathbb{R}^d$, $(A_s)_{s \in [t,T]}$ is an
$\mathbb{R}^d$-valued, $\mathbb{F}^{t,T}$-progressively measurable process,
$(\Gamma_s)_{s \in [t,T]}$ is an
${\cal S}^d$-valued, $\mathbb{F}^{t,T}$-progressively
measurable process such that
\be \label{Zbound}
\max \crl{|Z_s| \, , \, |A_s| \, , \, |\Gamma_s|}
\le m(1 + |X^{t,x}_s|^{p_4}) \quad \mbox{for all } s \in [t,T] \, ,
\ee
and
\be \label{gammalip}
|\Gamma_r - \Gamma_s|
\le m (1 + |X^{t,x}_r|^{p_5} + |X^{t,x}_s|^{p_5}) (|r-s| + |X^{t,x}_r - X^{t,x}_s|)
\quad \mbox{for all } r,s \in [t,T] \, .
\ee
Set ${\cal A}^{t,x} := \bigcup_{m \ge 0} {\cal A}_m^{t,x}$. It follows from
the assumptions (A1) and (A2) on $f$ and condition \eqref{Zbound} on $Z$ that for
all $y \in \mathbb{R}$ and $Z \in {\cal A}^{t,x}$, the forward SDE
 \bea \label{sdeY}
 d Y_s &=& f(s, X^{t,x}_s, Y_s, Z_s,\Gamma_s) \,ds +
           Z_s' \circ dX^{t,x}_s \, ,
 \quad s \in [t,T] \, ,\\
 Y_t &=& y \, ,
 \eea
has a unique strong solution $\left(Y^{t,x,y,Z}_s\right)_{s \in
[t,T]}$ (this can, for instance, be shown with the arguments in
the proofs of Theorems 2.3, 2.4 and 3.1 in Chapter IV of Ikeda
and Watanabe, 1989).

\subsection{Auxiliary stochastic target problems}

For every $m \ge 0$, we define the functions $V^m, U_m : [0,T]
\times \mathbb{R}^d \to \mathbb{R}$ as follows:
 \beas
 V^m(t,x) &:=& \inf \{y \in \mathbb{R} : Y^{t,x,y,Z}_T \ge g(X^{t,x}_T)
 \mbox{ for some } Z \in {\cal A}^{t,x}_m \} \, ,
 \eeas
 and
 \beas
 U_m(t,x) &:=& \sup \{y \in \mathbb{R} : Y^{t,x,y,Z}_T \le g(X^{t,x}_T)
 \mbox{ for some } Z \in {\cal A}_m^{t,x}\} \, .
 \eeas
Notice that these problems do not fit into the class of stochastic
target problems studied by Soner and Touzi (2002) and are more in
the spirit of Cheridito et al. (2005a, b).

 \begin{lemma} \label{pdpp} {\rm (Partial Dynamic Programming Principle)}\\
 Let $(t,x,m) \in [0,T) \times \mathbb{R}^d \times \mathbb{R}_+$ and
 $(y,Z) \in \mathbb{R} \times {\cal A}_m^{t,x}$ such that
 $Y^{t,x,y,Z}_T \ge g(X^{t,x}_T)$. Then
 \beas
 Y^{t,x,y,Z}_s
 &\ge&
 V^m(s, X^{t,x}_s) \quad \mbox{for all } s \in (t,T) \, .
 \eeas
\end{lemma}
\begin{proof}
Fix $s \in (t,T)$ and denote by $C^d[t,s]$ the set of all continuous functions
from $[t,s]$ to $\mathbb{R}^d$. Since $X^{t,x}_s$ and
$Y_s^{t,x,y,Z}$ are ${\cal F}^t_s$-measurable, there exist measurable functions
$$
\xi : C^d[t,s] \to \mathbb{R}^d \quad \mbox{and} \quad
\psi : C^d[t,s] \to \mathbb{R}
$$
such that
$$
X^{t,x}_s = \xi(W^{t,s}) \quad \mbox{and} \quad Y_s^{t,x,y,Z} = \psi(W^{t,s}) \, ,
$$
where we denote $W^{t,s} := (W^t_r)_{r \in [t,s]}$.
The process $Z$ is of the form
$$
Z_r = z + \int_t^r A_u du + \int_t^r \Gamma_u d X^{t,x}_u \, , \quad r \in [t,T] \, ,
$$
for $z \in \mathbb{R}$, $(A_r)_{r \in [t,T]}$ an
$\mathbb{R}^d$-valued, $\mathbb{F}^{t,T}$-progressively measurable process, and
$(\Gamma_r)_{r \in [t,T]}$ an ${\cal S}^d$-valued, $\mathbb{F}^{t,T}$-progressively
measurable process. Therefore, there exist progressively measurable functions
(see Definition 3.5.15 in Karatzas and Shreve, 1991)
\beas
\zeta, \phi &:& [t,T] \times C^d[t,T] \to \mathbb{R}^d\\
\chi &:& [t,T] \times C^d[t,T] \to {\cal S}^d
\eeas
such that
$$
Z_r = \zeta(r,W^{t,T}) \, , \quad
A_r = \phi(r,W^{t,T}) \quad \mbox{and} \quad
\Gamma_r = \chi(r,W^{t,T}) \quad \mbox{for } r \in [t,T] \, ,
$$
where $W^{t,T} := (W^t_r)_{r \in [t,T]}$.
With obvious notation, we define for every $w \in C^d[t,s]$ the $\mathbb{R}^d$-valued,
$\mathbb{F}^{s,T}$-progressively measurable process $(Z^w_r)_{r \in [s,T]}$ by
$$
Z^w_r = \zeta(s,w) + \int_s^r \phi(u, w + W^{s,T}) du +
\int_s^r \chi(u , w + W^{s,T}) du \, , \quad r \in [s,T] \, .
$$
Let $\mu$ be the distribution of $W^{t,s}$ on $C^d[t,s]$.
Then, $Z^w \in {\cal A}^{s, \xi(w)}_m$ for $\mu$-almost all $w \in C^d[t,s]$, and
\beas
1 &=& P \edg{Y^{t,x,y,Z}_T \ge g(X^{t,x}_T)}\\ &=&
\int_{C^d[t,s]} P \edg{Y^{s, \xi(w), \psi(w), Z^w}_T \ge g(X^{s, \xi(w)}_T)
\mid W^{t,s} = w} d \mu(w) \, .
\eeas
Hence, for $\mu$-almost all $w \in C^d[t,s]$, the control $Z^w$ satisfies
$$
P \edg{Y_T^{s, \xi(w) , \psi(w) , Z^w} \ge g(X_T^{s, \xi(w)})} =
1 \,,
$$
which shows that $\psi(w)=Y_s^{s,\xi(w),\psi(w),Z^w}\ge
V^m\left(s, \xi(w)\right)$. In view of the definition of the
functions $\xi$ and $\psi$, this implies that $Y^{t,x,y,Z}_s \ge
V^m(s, X^{t,x}_s)$.
 \end{proof}

Since we have no a priori knowledge of any regularity of the
functions $V^m$ and $U_m$, we introduce the semi-continuous
envelopes as in Barles and Perthame (1988)
 \beas
 V^m_*(t,x)
 \;:=\;
 \liminf_{(\tilde t,\tilde x)\to(t,x)}\;
 V^m\left(\tilde t,\tilde x\right)
 &\mbox{and}&
 U_m^*(t,x)
 \;:=\;
 \limsup_{(\tilde t,\tilde x)\to(t,x)}\;
 U_m\left(\tilde t,\tilde x\right) \,,
 \eeas
for all $(t,x)\in[0,T]\times\mathbb{R}^d$, where we set by
convention $V^m(t,x)=-U^m(t,x)=\infty$ for $t \notin [0,T]$. For the theory
of viscosity solutions, we refer to Crandal et al. (1992) and the
book of Fleming and Soner (1993).

 \begin{theorem} \label{supersolution}
Let $m \ge 0$, and assume that $V^m_*$ is $\mathbb{R}$-valued.
Then $V^m_*$ is a viscosity supersolution of the {\rm PDE}
 \beas
 - v_t(t,x)
 \,+\,
 \sup_{\beta \in {\cal S}^d_+}\; f\left(t,x,v(t,x),Dv(t,x),D^2v(t,x)+\beta\right)
 &=& 0
 \quad \mbox{on } [0,T) \times \mathbb{R}^d \, .
 \eeas
 \end{theorem}

Before turning to the proof of this result, let us state the
corresponding claim for the value function $U_m$.

\begin{corollary} \label{subsolution}
Let $m \ge 0$, and assume that $U_m^*$ is
$\mathbb{R}$-valued. Then $U_m^*$ is a viscosity subsolution of the {\rm
PDE}
 \bea \label{infpde}
 - u_t(t,x)
 \,+\,
 \inf_{\beta \in {\cal S}^d_+}\; f\left(t,x,u(t,x),Du(t,x),D^2u(t,x)-\beta\right)
 &=& 0
 \quad \mbox{on } [0,T) \times \mathbb{R}^d \, .
 \eea
\end{corollary}

\begin{proof}
 Observe that for all $(t,x) \in [0,T) \times \mathbb{R}^d$,
 \beas
 - U_m(t,x)
 &=&
 \inf \left\{y \in \mathbb{R} : \hat{Y}^{t,x,y,Z}_T \ge -
 g(X^{t,x}_T) \mbox{ for some } Z \in {\cal A}_m^{t,x} \right\} \, ,
 \eeas
where for given $(y,Z) \in \mathbb{R} \times {\cal A}_m^{t,x}$,
the process $\hat{Y}^{t,x,y,Z}$ is the unique strong solution of
the SDE
 $$
 Y_s = y + \int_t^s - f(r,X^{t,x}_r, - Y_r, -Z_r,
 -\Gamma_r) dr + \int_t^s (Z_r)' \circ dX^{t,x}_r \, , \quad s \in [t,T]
 \, .
 $$
Hence, it follows from Theorem \ref{supersolution} that $- U_m^*
= (- U_m)_*$ is a viscosity supersolution of the PDE
 \beas
 - u_t(t,x) \,-\, \inf_{\beta \in {\cal S}^d_+}\;
 f\left(t, x, -u(t,x), - Du(t,x), - D^2u(t,x) - \beta\right)
 &=& 0
 \quad \mbox{on } [0,T) \times \mathbb{R}^d \, ,
 \eeas
which shows that $U_m^*$ is a viscosity subsolution of the PDE
\eqref{infpde} on $[0,T) \times \mathbb{R}^d$.
\end{proof}

\vspace{5mm}

\noindent {\bf Proof of Theorem \ref{supersolution}}
Choose $(t_0,x_0) \in [0,T) \times \mathbb{R}^d$ and
 $\varphi \in C^{\infty}([0,T) \times \mathbb{R}^d)$ such
 that
 $$
 0 = (V^m_* - \varphi) (t_0,x_0) = \min_{(t,x) \in [0,T)
 \times \mathbb{R}^d} (V^m_* - \varphi)(t,x) \, .
 $$
 Let $(t_n , x_n) \to (t_0,x_0)$ such that $V^m(t_n,x_n) \to
 V^m_*(t_0,x_0)$. There exist positive numbers
 $\varepsilon_n \to 0$ such that for $y_n = V^m(t_n,x_n) + \varepsilon_n$,
 there exists $Z^n \in {\cal A}_m^{t_n,x_n}$ with
 \beas
 Y^n_T &\ge& g(X^n_T) \, ,
 \eeas
 where we denote $(X^n, Y^n) = (X^{t_n, x_n}, Y^{t_n, x_n,y_n,Z^n})$ and
 \beas
 Z^n_s &=& z_n + \int_{t_n}^s A^n_r dr + \int_{t_n}^s \Gamma^n_r dX^n_r
 \, , \quad s \in [t_n,T] \, .
 \eeas
 Note that for all $n$, $\Gamma^n_{t_n}$ is almost surely constant, and
 $|z_n|, |\Gamma^n_{t_n}| \le m(1 + |x_n|^{p_4})$ by assumption \eqref{Zbound}.
 Hence, by passing
 to a subsequence, we can assume that $z_n \to z_0 \in
 \mathbb{R}^d$ and $\Gamma^n_{t_n} \to \gamma_0 \in {\cal S}^d$.
 Observe that $\alpha_n := y_n - \varphi(t_n,x_n) \to 0$. We choose
 a decreasing sequence of numbers $\delta_n \in (0, T - t_n)$ such that
 $\delta_n \to 0$ and $\alpha_n / \delta_n \to 0$. By Lemma \ref{pdpp},
 \beas
 Y^n_{t_n + \delta_n}
 &\ge&
 V^m\left(t_n + \delta_n , X^n_{t_n + \delta_n}\right) \, ,
 \eeas
 and therefore,
 \beas
 Y^n_{t_n + \delta_n} - y_n + \alpha_n
 &\ge&
 \varphi\left(t_n + \delta_n,X^n_{t_n + \delta_n}\right)
 - \varphi(t_n,x_n) \, ,
 \eeas
 which, after two applications of It\^{o}'s Lemma, becomes
 \bea
 \nonumber
 \alpha_n &+& \int_{t_n}^{t_n + \delta_n}
 [f(s, X^n_s, Y^n_s, Z^n_s, \Gamma^n_s) - \varphi_t(s,X^n_s)] ds\\
 \nonumber
 &+& [z_n - D \varphi(t_n,x_n)]' [X^n_{t_n + \delta_n} - x_n]\\
 \nonumber
 &+& \int_{t_n}^{t_n + \delta_n} \left( \int_{t_n}^s A^n_r -
 {\cal L} D \varphi(r, X^n_r) dr \right)' \circ dX^n_s\\
 \label{inequality}
 &+& \int_{t_n}^{t_n + \delta_n} \left( \int_{t_n}^s \Gamma^n_r -
 D^2 \varphi(r, X^n_r) dX^n_r \right)' \circ dX^n_s
 \;\ge\; 0
 \eea
 It is shown in Lemma \ref{convergence} below that the sequence of random vectors
 \be \label{vectn}
 \left(
 \begin{array}{c}
 \delta_n^{-1} \int_{t_n}^{t_n + \delta_n}
 [f(s, X^n_s, Y^n_s, Z^n_s, \Gamma^n_s) - \varphi_t(s, X^n_s)] ds\\
 \delta_n^{-1/2} [X^n_{t_n + \delta_n} - x_n]\\
 \delta_n^{-1} \int_{t_n}^{t_n+\delta_n} \left( \int_{t_n}^s A^n_r
 - {\cal L} D \varphi(r, X^n_r) dr \right)' \circ dX^n_s \\
 \delta_n^{-1} \int_{t_n}^{t_n + \delta_n} \left( \int_{t_n}^s \Gamma^n_r
 - D^2 \varphi(r , X^n_r) dX^n_r \right)' \circ dX^n_s
 \end{array}
 \right) \, , \, n \ge 1 \, ,
 \ee
 converges in distribution to
 \be \label{vect}
 \left(
 \begin{array}{c}
 f(t_0, x_0, \varphi(t_0,x_0) , z_0, \gamma_0) - \varphi_t(t_0,x_0)\\
 \sigma(x_0) W_1\\
 0\\
 \frac{1}{2} W_1' \sigma(x_0)' [\gamma_0 - D^2 \varphi(t_0,x_0)]
 \sigma(x_0) W_1
 \end{array}
 \right)
 \ee
 Set $\eta_n = |z_n - D \varphi(t_n,x_n)|$, and assume
 $\delta_n^{-1/2} \eta_n \to \infty$ along a subsequence.
 Then, along another subsequence, $\eta_n^{-1} (z_n - D
 \varphi(t_n,x_n))$ converges to some $\eta_0 \in \mathbb{R}^d$ with
 \be \label{1}
 |\eta_0| =1 \, .
 \ee
 Multiplying inequality \eqref{inequality} with
 $\delta_n^{-1/2} \eta_n^{-1}$ and passing to the limit yields
 $$
 \eta_0' \sigma(x_0) W_1 \ge 0 \, ,
 $$
 which, since $\sigma(x_0)$ is invertible, contradicts \eqref{1}.
 Hence, the sequence $(\delta_n^{-1/2} \eta_n)$ has to be bounded,
 and therefore, possibly after passing to a subsequence,
 $$
 \delta_n^{-1/2} [z_n - D \varphi(t_n,x_n)]
 \quad \mbox{converges to some } \xi_0 \in \mathbb{R}^d \, .
 $$
 It follows that $z_0 = D \varphi(t_0,x_0)$. Moreover, we can
 divide inequality \eqref{inequality} by $\delta_n$ and pass to the
 limit to get
 \be \label{W1}
 \begin{aligned}
 & f(t_0, x_0, \varphi(t_0,x_0) , D \varphi(t_0,x_0) , \gamma_0) - \varphi_t(t_0,x_0)\\
 + \; & \xi_0' \sigma(x_0) W_1 + \frac{1}{2} W_1' \sigma(x_0)'
 [\gamma_0 - D^2 \varphi(t_0,x_0)]
 \sigma(x_0) W_1
 \;\ge\; 0 \, .
 \end{aligned}
 \ee
 Since the support of the random vector $W_1$ is $\mathbb{R}^d$, it follows from
 \eqref{W1} that
 \beas
 && f(t_0,x_0, \varphi(t_0,x_0) , D \varphi(t_0,x_0), \gamma_0) - \varphi_t(t_0,x_0)\\
 &+& \xi_0' \sigma(x_0) w + \frac{1}{2} w' \sigma(x_0)'
 [\gamma_0 - D^2 \varphi(t_0,x_0)]
 \sigma(x_0) w \ge 0 \, ,
 \eeas
for all $w \in \mathbb{R}^d$. This shows that
 $$
 f(t_0,x_0, \varphi(t_0,x_0) , D \varphi(t_0,x_0), \gamma_0) -
 \varphi_t(t_0,x_0) \ge 0
 \quad \mbox{and} \quad
 \beta := \gamma_0-D^2 \varphi(t_0,x_0)\ge 0 \, ,
$$
and hence,
$$
- \varphi_t(t_0,x_0) + \sup_{\beta \in {\cal S}^d_+}
f(t_0,x_0, \varphi(t_0,x_0) , D \varphi(t_0,x_0), D^2 \varphi(t_0,x_0) + \beta) \ge 0 \, .
$$
 \ep

\begin{lemma} \label{convergence}
The sequence of random vectors \eqref{vectn} converges in distribution to \eqref{vect}.
\end{lemma}
\begin{proof}
With the methods used to solve Problem 5.3.15 in Karatzas and
Shreve (1991) it can also be shown that for all fixed $q
> 0$ and $m \ge 0$, there exists a constant $C \ge 0$ such that
for all $0 \le t \le s \le T$, $x \in \mathbb{R}^d$, $y \in
\mathbb{R}$ and $Z \in {\cal A}^{t,x}_m$, \bea
\label{Xuniformbound}
\E{\max_{s \in [t,T]} |X^{t,x}_s|^q} &\le& C(1 + |x|^q)\\
\label{Xincrements}
\E{\max_{r \in [t,s]} |X^{t,x}_r - x|^q} &\le& C(1 + |x|^q) (s-t)^{q/2}\\
\label{Yuniformbound}
\E{\max_{s \in [t,T]} |Y^{t,x,y,Z}_s|^q} &\le& C(1 + |y|^q + |x|^{\tilde{q}})\\
\label{Yincrements}
\E{\max_{r \in [t,s]} |Y^{t,x,y,Z}_r - y|^q} &\le& C(1 + |y|^q
+ |x|^{\tilde{q}}) (s-t)^{q/2} \, ,
\eea
where $\tilde{q} := \max \crl{p_2 q \, , \, p_2 p_4 q \, , \, (p_4 + 2 p_1) q}$.
For every $n \ge 1$, we introduce the $\mathbb{F}^{t_n,T}$-stopping time
$$
\tau_n := \inf \{s \ge t_n : X^n_s \notin B_1(x_0)\} \wedge (t_n + \delta_n) \, ,
$$
where $B_1(x_0)$ denotes the open unit ball
in $\mathbb{R}^d$ around $x_0$. It follows from the fact that
$x_n \to x_0$ and \eqref{Xincrements} that
\be \label{taunconv}
P [\tau_n < t_n + \delta_n] \to 0 \, .
\ee
The difference
$$
(X^n_{t_n + \delta_n} - x_n) - \sigma(x_0) (W_{t_n + \delta_n} - W_{t_n})
$$
can be written as
$$
\int_{t_n}^{t_n + \delta_n} \mu(X^n_s) ds + \int_{t_n}^{t_n + \delta_n}
[\sigma(X^n_s) - \sigma(x_n)] dW_s + (\sigma(x_n) - \sigma(x_0)) (W_{t_n + \delta_n} - W_{t_n}) \, ,
$$
and obviously,
$$
\frac{1}{\sqrt{\delta_n}} (\sigma(x_n) - \sigma(x_0)) (W_{t_n + \delta_n} - W_{t_n}) \to 0
\quad \mbox{in } L^2 \, .
$$
Moreover, it can be deduced with standard arguments from \eqref{lip}, \eqref{growth},
\eqref{Xuniformbound} and \eqref{Xincrements} that
$$
\frac{1}{\sqrt{\delta_n}} \int_{t_n}^{t_n + \delta_n} \mu(X^n_s) ds \to 0 \quad
\mbox{and} \quad
\frac{1}{\sqrt{\delta_n}} \int_{t_n}^{\tau_n}
[\sigma(X^n_s) - \sigma(x_n)] dW_s \to 0 \quad \mbox{in } L^2 \, .
$$
This shows that
\be \label{term1}
\frac{1}{\sqrt{\delta_n}}
\crl{X^n_{t_n + \delta_n} - x_n - \sigma(x_0) (W_{t_n + \delta_n} - W_{t_n})}
\to 0 \quad \mbox{in probability}.
\ee
Similarly, it can be derived from the boundedness assumption \eqref{Zbound} on $A^n$
that
\be \label{term2}
\frac{1}{\delta} \int_{t_n}^{t_n + \delta_n} \left( \int_{t_n}^s A^n_r -
{\cal L} D \varphi(r,X^n_r) \, dr \right)' \circ
dX^n_s \to 0 \quad \mbox{in probability.}
\ee
By the continuity assumption \eqref{gammalip} on $\Gamma^n$,
$$
\frac{1}{\delta_n} \int_{t_n}^{t_n + \delta_n} \left( \int_{t_n}^s [ \Gamma^n_r
- \Gamma^n_{t_n}] dX^n_r \right)' \circ dX^n_s \to 0 \quad \mbox{in } L^2 \, ,
$$
and
\beas
\frac{1}{\delta_n} && \hspace*{-5mm} \left\{ \int_{t_n}^{t_n + \delta_n}
\left( \int_{t_n}^s \Gamma^n_{t_n} dX^n_r \right)' \circ dX^n_s \right.\\
&& \left. - \frac{1}{2} (W_{t_n + \delta_n} - W_{t_n})' \sigma(x_0)' \gamma_0
\sigma(x_0) (W_{t_n + \delta_n} - W_{t_n}) \right\} \to 0
\quad \mbox{in probability}.
\eeas
Hence,
\bea
\label{term3}
\frac{1}{\delta_n} && \hspace*{-5mm} \left\{ \int_{t_n}^{t_n + \delta_n}
\left( \int_{t_n}^s \Gamma^n_r dX^n_r \right)'
\circ dX^n_s \right. \\
\nonumber
&& \left. - \frac{1}{2} (W_{t_n + \delta_n} - W_{t_n})' \sigma(x_0)' \gamma_0 \sigma(x_0)
(W_{t_n + \delta_n} - W_{t_n}) \right\} \to 0 \quad \mbox{in probability}.
\eea
Similarly, it can be shown that
\bea
\label{term4}
\frac{1}{\delta_n} && \hspace*{-5mm} \left\{
\int_{t_n}^{t_n + \delta_n} \left( \int_{t_n}^s D^2
\varphi(r, X^n_r) dX^n_r \right)' \circ dX^n_s \right.\\
\nonumber
&& \left. - \frac{1}{2} (W_{t_n + \delta_n} - W_{t_n})' \sigma(x_0)'
D^2 \varphi(t_0,x_0) \sigma(x_0) (W_{\tau_n} - W_{t_n}) \right\}
\to 0 \quad \mbox{in probability}.
\eea
Finally, it follows from the continuity of $f$ and $\varphi_t$ as well as
\eqref{Zbound}, \eqref{gammalip},
\eqref{Xuniformbound}, \eqref{Xincrements} and \eqref{Yincrements} that
\be \label{term5}
\frac{1}{\delta_n} \int_{t_n}^{t_n + \delta_n}
[f(X^n_s, Y^n_s, Z^n_s, \Gamma^n_s) - \varphi_t(s, X^n_s)] \, ds
\to f(x_0, D \varphi(t_0,x_0) , z_0, \gamma_0) - \varphi_t(t_0, x_0)
\ee
in probability.
Now, the lemma follows from \eqref{term1}--\eqref{term5} and the simple fact
 that for each $n$, the random vector
 $$
 \left(
 \begin{array}{c}
 f(x_0, \varphi(t_0,x_0) , z_0, \gamma_0) - \varphi_t(t_0,x_0)\\
 \delta_n^{-1/2} \sigma(x_0) (W_{t_n + \delta_n} - W_{t_n})\\
 0\\
 \delta_n^{-1} \frac{1}{2} (W_{t_n + \delta_n} - W_{t_n})' \sigma(x_0)'
 [\gamma_0 - D^2 \varphi(t_0,x_0)] \sigma(x_0) (W_{t_n + \delta_n} - W_{t_n})
 \end{array}
 \right)
 $$
 has the same distribution as
 $$
 \left(
 \begin{array}{c}
f(x_0, \varphi(t_0,x_0) , z_0, \gamma_0) - \varphi_t(t_0,x_0)\\
\sigma(x_0) W_1\\
0\\
 \frac{1}{2} W_1' \sigma(x_0)' [\gamma_0 - D^2 \varphi(t_0,x_0)]
 \sigma(x_0) W_1
 \end{array}
 \right)
 $$
 \end{proof}

We conclude this subsection by the following bounds on the growth of
the value functions $V^m$ and $U_m$.

\begin{proposition} \label{Vbound}
Let $p = \max \crl{p_2 \, , \, p_3 \, , \, p_2 p_4 \, , \, p_4 + 2 p_1}$. Then
there exists for every $m \ge 0$ a constant $C_m \ge 0$ such that
\bea
\label{lowerbound}
V^m_*(t,x) &\ge& - C_m(1 + |x|^p) \quad \mbox{and}\\
\label{upperbound}
U_m^*(t,x) &\le& C_m(1 + |x|^p)
\eea
for all $(t,x) \in [0,T] \times \mathbb{R}^d$. Moreover,
\bea
\label{lowerboundT}
V^m_*(T,x) &\ge& g(x) \quad \mbox{and}\\
\label{upperboundT}
U_m^*(T,x) &\le& g(x)
\eea
for all $x \in \mathbb{R}$.
\end{proposition}

\begin{proof}
We show \eqref{upperbound} and \eqref{upperboundT}. The proofs of
\eqref{lowerbound} and \eqref{lowerboundT} are completely
analogous. To prove \eqref{upperbound} it is enough to show that
for fixed $m \ge 0$,
there exists a constant $C_m \ge 0$ such that for all $(t,x) \in
[0,T] \times \mathbb{R}^d$ and $(y,Z) \in \mathbb{R} \times {\cal
A}^{t,x}_m$ satisfying $Y^{t,x,y,Z}_T \le g(X^{t,x}_T)$, we have
 \beas
 y &\le& C_m(1+|x|^p) \,.
 \eeas
For $y \le 0$ there is nothing to show. So, we assume $y > 0$ and
introduce the stopping time
$$
\tau := \inf \crl{s \ge t \mid Y^{t,x,y,Z}_s = 0} \wedge T \, .
$$
Then, we have for all $s \in [t,T]$,
\beas
&& Y^{t,x,y,Z}_{s \wedge \tau} + \int_s^{s \vee \tau}
f(r,X^{t,x}_r, Y^{t,x,y,Z}_r, Z_r, \Gamma_r) dr
+ \int_s^{s \vee \tau} Z_r' \mu(X^{t,x}_r) dr\\
&&+ \int_s^{s \vee \tau} Z_r' \sigma(X^{t,x}_r) dW_r
+ \frac{1}{2} \int_s^{s \vee \tau} {\rm Tr} [\Gamma_r
\sigma(X^{t,x}_r) \sigma(X^{t,x}_r)'] dr\\
&=& Y^{t,x,y,Z}_{s \wedge \tau} + \int_s^{s \vee \tau}
f(r,X^{t,x}_r, Y^{t,x,y,Z}_r, Z_r, \Gamma_r) dr
+ \int_s^{s \vee \tau} Z_r' \circ d X^{t,x}_r\\
&=& Y^{t,x,y,Z}_{\tau}\\
&\le& g(X^{t,x}_T) \vee 0 \, .
\eeas
Hence, it follows from (A2), (A3), \eqref{growth} and \eqref{Zbound} that
for $\tilde{p} = \max \crl{p_2 \, , \, p_2 p_4 \, , \, p_4 + 2 p_1}$,
\beas
&& h(s) := \E{Y^{t,x,y,Z}_{s \wedge \tau}}\\
&\le& \E{\abs{g(X^{t,x}_T)}} + \E{\int_s^{s \vee \tau}
\abs{f(r,X^{t,x}_r, Y^{t,x,y,Z}_r, Z_r, \Gamma_r)} dr}\\
&& + \E{\int_s^{s \vee \tau} \abs{Z_r' \mu(X^{t,x}_r)} dr}
+ \E{\frac{1}{2} \int_s^{s \vee \tau} \abs{{\rm Tr} [\Gamma_r
\sigma(X^{t,x}_r) \sigma(X^{t,x}_r)']} dr}\\
&\le& G \, \E{1 + |X^{t,x}_T|^{p_3}} + F \int_s^T h(r) dr
+ K \int_s^T (1 + |X^{t,x}_r|^{\tilde{p}}) dr\\
&\le& \tilde{K} (1 + \abs{x}^p) + F \int_s^T h(r) dr \, ,
\eeas
for constants $K$ and $\tilde{K}$ independent of $t$, $x$, $y$ and $Z$.
It follows from Gronwall's lemma that
$$
h(s) \le \tilde{K} (1 + \abs{x}^p) e^{F(T-s)} \quad \mbox{for all }
s \in [t,T] \, .
$$
In particular, $y = h(t) \le C_m(1 + |x|^p)$ for some constant $C_m$ independent of
$t$, $x$, $y$ and $Z$.

To prove \eqref{upperboundT}, we assume by way of contradiction
that there exists an $x \in \mathbb{R}^d$ such that $U^*_m(T,x)
\ge g(x) + 3 \varepsilon$ for some $\varepsilon > 0$. Then, there
exists a sequence $(t_n, x_n)_{n \ge 1}$ in $[0,T) \times
\mathbb{R}^d$ converging to $(T,x)$ such that $U_m(t_n, x_n) \ge
g(x) + 2 \varepsilon$ for all $n \ge 1$. Hence, for every integer
$n \ge 1$, there exists a real number $y_n \in [g(x) +
\varepsilon, g(x) + 2 \varepsilon]$ and a process $Z^n \in {\cal
A}^{t_n,x_n}_m$ of the form $Z^n_s = z_n + \int_{t_n}^s A^n_s ds
+ \int_{t_n}^s \Gamma^n_s d X^{t_n,x_n}_s$ such that
 \be
 \label{upperboundonyn}
 y_n
 \le
 g(X^{t_n,x_n}_T)
 - \int_{t_n}^T f(s,X^{t_n,x_n}_s, Y^{t,x_n,y_n,Z^n}_s, Z^n_s, \Gamma^n_s) ds
 - \int_{t_n}^T (Z^n_s)' \circ d X^{t_n,x_n}_s \, .
 \ee
By \eqref{Zbound}, \eqref{Xuniformbound}, \eqref{Xincrements} and
\eqref{Yuniformbound}, the right-hand side of
\eqref{upperboundonyn} converges to $g(x)$ in probability.
Therefore, it follows from \eqref{upperboundonyn} that $g(x) +
\varepsilon \le g(x)$. But this is absurd, and hence, we must
have $U^*_m(T,x) \le g(x)$ for all $x \in \mathbb{R}^d$.
\end{proof}

 \subsection{Main result}

For our main result, Theorem \ref{thmmain} below, we need two more
assumptions on the functions $f$ and $g$, the first of which is

{\bf (A4)} \quad For all $(t,x,y,z) \in [0,T] \times \mathbb{R}^d
\times \mathbb{R} \times \mathbb{R}^d$ and $\gamma, \tilde{\gamma}
\in {\cal S}^d$,
 \beas
 f(t,x,y,z,\gamma) \;\ge\; f(t,x,y,z, \tilde{\gamma})
 &\mbox{whenever}&
 \gamma \;\le\;\tilde{\gamma} \,.
 \eeas

\begin{Remark} \label{viscsol}
Under (A1)--(A4) it immediately follows from Theorem \ref{supersolution} that
$V^m_*$ is a viscosity supersolution of the PDE \eqref{pde} on $[0,T) \times \mathbb{R}^d$,
provided it is $\mathbb{R}$-valued. Analogously, if (A1)--(A4) hold and
$U_m^*$ is $\mathbb{R}$-valued, Corollary \ref{subsolution} implies that $U_m^*$ is a
viscosity subsolution of the PDE \eqref{pde} on $[0,T) \times \mathbb{R}^d$.
\end{Remark}

For our last assumption and the statement of Theorem \ref{thmmain}
we need the following

\begin{definition}
Let $q \ge 0$.\\
{\bf 1.} We call a function $v : [0,T] \times \mathbb{R}^d \to \mathbb{R}$ a viscosity solution
with growth $q$ of the {\rm PDE} \eqref{pde} with terminal condition \eqref{terminal}
if $v$ is a viscosity solution of \eqref{pde} on $[0,T) \times \mathbb{R}^d$ such that
$v^*(T,x) = v_*(T,x) = g(x)$ for all $x \in \mathbb{R}^d$ and there exists a
constant $C$ such that $\abs{v(t,x)} \le C(1 + \abs{x}^q)$
for all $(t,x) \in [0,T] \times \mathbb{R}^d$.\\
{\bf 2.}
We say that the {\rm PDE} \eqref{pde} with terminal condition \eqref{terminal}
has comparison with growth $q$ if the following holds:

If $w : [0,T] \times \mathbb{R}^d \to \mathbb{R}$ is lower semicontinuous and a
viscosity supersolution of \eqref{pde} on $[0,T) \times \mathbb{R}^d$ and
$u : [0,T] \times \mathbb{R}^d \to \mathbb{R}$ upper semicontinuous and a
viscosity subsolution of \eqref{pde} on $[0,T) \times \mathbb{R}^d$ such that
$$
w(T,x) \ge g(x) \ge u(T,x) \quad \mbox{for all } x \in
\mathbb{R}^d
$$
and there exists a constant $C \ge 0$ with
$$
w(t,x) \ge - C (1 + |x|^p) \quad \mbox{and}
\quad u(t,x) \le C (1 + |x|^p) \quad \mbox{for all }
(t,x) \in [0,T) \times \mathbb{R}^d \, ,
$$
then
$w \ge u$ on $[0,T] \times \mathbb{R}^d$.
\end{definition}
With this definition our last assumption on $f$ and $g$ is

{\bf (A5)} \quad The PDE \eqref{pde} with terminal condition
\eqref{terminal} has comparison with growth $p = \max \crl{p_2 \,
, \, p_3 \, , \, p_2 p_4 \, , \, p_4 + 2 p_1}$.

\begin{Remarks} $\mbox{}$\\
\noindent {\bf 1.} The monotonicity assumption (A4) is natural
from the PDE viewpoint. It implies that $f$ is elliptic and the
PDE \eqref{pde} parabolic. If $f$ satisfies the following stronger
version of (A4): there exists a constant $C > 0$ such that \be
\label{uniformelliptic} f(t,x,y,z, \gamma - B) \ge
f(t,x,y,z,\gamma) + C \ {\rm Tr}[B] \ee for all $(t,x,y,z,\gamma)
\in [0,T) \times \mathbb{R}^d \times \mathbb{R} \times
\mathbb{R}^d \times {\cal S}^d$ and $B \in {\cal S}^d_+$, then the
PDE \eqref{pde} is called uniformly parabolic, and there exist
general results on existence, uniqueness and smoothness of
solutions, see for instance Krylov (1987) or Evans (1998). When
$f$ is linear in the $\gamma$ variable (in particular, for the
semi- and quasi-linear equations discussed in Subsections
\ref{sectnumericalsl} and \ref{sectnumericalql} below), the
condition \eqref{uniformelliptic} essentially guarantees
existence, uniqueness and smoothness of solutions to the PDE
\eqref{pde}--\eqref{terminal}; see for instance, Section 5.4 in Ladyzenskaya et al.
(1967). Without parabolicity there are no comparison results for
PDEs of the form \eqref{pde}--\eqref{terminal}.
\\
{\bf 2.} Condition (A5) is an implicit assumption on the functions
$f$ and $g$. But we find it more convenient to assume comparison
directly in the form (A5) instead of putting technical assumptions
on $f$ and $g$ which guarantee that the PDE \eqref{pde} with
terminal condition \eqref{terminal} has comparison. Several
comparison results for non-linear PDEs are available in the
literature; see for example, Crandall et al. (1992), Fleming and
Soner (1993), Cabre and Caffarelli (1995). However, most results
are stated for equations in bounded domains. For equations in the
whole space, the critical issue is the interplay between the
growth of solutions at infinity and the growth of the
non-linearity. We list some typical situations where comparison
holds:

{\bf a)} {\sl Comparison with growth 1:} Assume (A1)--(A4) and there exists a function
$h :[0,\infty] \to [0,\infty]$ with $\lim_{x \to 0} h(x) = 0$ such that
$$
\left| f(t,x,y,\alpha (x-\tilde x),A) - f(t,\tilde x,y,\alpha
(x-\tilde x),B) \right| \le h( \alpha |x-\tilde x|^2 +
|x-\tilde x|),
$$
for all $(t,x,\tilde x, y)$, $\alpha> 0$ and $A$, $B$ satisfying
$$
-\alpha \left[
\begin{array}{rr}
I & 0 \\
0 & I
\end{array}
\right]
\le
 \left[
\begin{array}{rr}
A & 0 \\
0 & -B
\end{array}
\right] \le
 \alpha \left[
\begin{array}{rr}
I & -I \\
-I & I
\end{array}
\right] .
$$
Then it follows from Theorem 8.2 in Crandall et al. (1992) that
equations of the form \eqref{pde}, \eqref{terminal} have
comparison with growth $0$ if the domain is bounded. If the domain
is unbounded, it follows from the modifications outlined in
Section 5.D of Crandall et al. (1992) that \eqref{pde} and
\eqref{terminal} have comparison with growth $1$.

{\bf b)} When the equation \eqref{pde} is a dynamic programming
equation related to a stochastic optimal control problem, then a
comparison theorem for bounded solutions is given in Fleming and
Soner (1993), Section 5.9, Theorem V.9.1. In this case, $f$ has
the form
$$
f(t,x,y,z,\gamma) = \sup_{u \in U}\ \crl{\alpha(t,x,u) + \beta(t,x,u) y
+ b(t,x,u)' z - {\rm Tr} [c(t,x,u) \gamma]} \, ,
$$
see Subsection \ref{sectnumericalstoch} below. Theorem V.9.1 in
Fleming and Soner (1993) is proved under the assumption that
$\beta \equiv 0$, $U$ is a compact set and that $\alpha, b, c$ are
uniformly Lipschitz and growing at most linearly (see IV (2.1) in
Fleming and Soner, 1993).  This result can be extended directly
to the case where $\beta$ satisfies a similar condition and to
equations related to differential games, that is, when
$$
f(t,x,y,z,\gamma) = \sup_{u \in U} \inf _{\tilde{u} \in \tilde{U}}
\left\{ \alpha(t,x,u,\tilde{u})
+ \beta(t,x,u, \tilde{u}) y + b(t,x,u,\tilde{u})' z - {\rm Tr}
[c(t,x,u,\tilde{u}) \gamma] \right\} \, .
$$

{\bf c)} Many techniques in dealing with unbounded solutions were
developed by Ishii (1984) for first order equations (that is, when
$f$  is independent of $\gamma$).  These techniques can be
extended to second order equations.  Some related
results can be found in Barles et al. (1997, 2003). In Barles et
al (1997), in addition to comparison results for PDEs, one can
also find BSDEs based on jump Markov processes.
\end{Remarks}

\begin{theorem}[Uniqueness of 2BSDE] \label{thmmain}
Assume {\rm (A1)--(A5)} and there exists an $x_0 \in \mathbb{R}^d$ such that
the {\rm 2BSDE} corresponding to $(X^{0,x_0},f,g)$ has a solution
$(Y^{0,x_0},Z^{0,x_0},\Gamma^{0,x_0},A^{0,x_0})$ with
$Z^{0,x_0} \in {\cal A}^{0,x_0}$. Then the following hold:\\
{\bf 1.}
The associated {\rm PDE} \eqref{pde} with terminal condition \eqref{terminal}
has a unique viscosity solution $v$ with growth
$p = \max \crl{p_2 \, , \, p_3 \, , \, p_2 p_4 \, , \, p_4 + 2 p_1}$,
and $v$ is continuous on $[0,T] \times \mathbb{R}^d$.\\
{\bf 2.} For all $(t,x) \in [0,T) \times \mathbb{R}^d$, there
exists exactly one solution $(Y^{t,x}, Z^{t,x}, \Gamma^{t,x},
A^{t,x})$ to the {\rm 2BSDE} corresponding to $(X^{t,x}, f, g)$
such that $Z^{t,x} \in {\cal A}^{t,x}$, and the process $Y^{t,x}$
satisfies
 \bea \label{representationY}
 Y^{t,x}_s &=& v(s, X^{t,x}_s) \, , \quad s \in [t,T] \,,
 \eea
where $v$ is the unique continuous viscosity solution with growth $p$ of
\eqref{pde}--\eqref{terminal}.
\end{theorem}

Before turning to the proof of this result, we make some remarks.

\begin{Remark}
If the assumptions of Theorem \ref{thmmain} are fulfilled,
it follows from \eqref{representationY} that
$v(t,x) = Y^{t,x}_t$
for all $(t,x) \in [0,T) \times \mathbb{R}^d$.
Hence, $v(t,x)$ can be approximated by backward simulation of the process
$(Y^{t,x}_s)_{s \in [t,T]}$.
If $v$ is $C^2$, it follows from It\^{o}'s lemma
that $Z^{t,x}_s = Dv(s, X^{t,x}_s)$, $s \in [t,T]$. Then,
$Dv(t,x)$ can be approximated by backward simulation of
$(Z^{t,x}_s)_{s \in [t,T]}$.
If $v$ is $C^3$, then also
$\Gamma^{t,x}_s = D^2 v(s, X^{t,x}_s)$, $s \in [t,T]$, and
$D^2v(t,x)$ can be approximated by backward simulation of
$(\Gamma^{t,x}_s)_{s \in [t,T]}$. A formal discussion of a
potential numerical scheme for the backward simulation of the processes $Y^{t,x}$,
$Z^{t,x}$ and $\Gamma^{t,x}$ is provided in Subsection
\ref{sectnumericalnl} below.
\end{Remark}

\begin{Remark}
Assume there exists a classical solution $v$ of the PDE \eqref{pde} such that
$$
v_t, Dv, D^2v, {\cal L} D v \mbox{ exist and are continuous on }
[0,T] \times \mathbb{R}^d \, ,
$$
and there exists a constant $m \ge 0$ such that
\be \label{vbound}
\left.
\begin{array}{l}
|Dv(t,x)|\\
|D^2v(t,x)|\\
|{\cal L} Dv (t,x)|
\end{array}
\right\}
\le m(1 + |x|^{p_4}) \quad \mbox{for all } 0 \le t \le T \mbox{
and } x \in \mathbb{R}^d \, .
\ee
and
\bea \label{vlip}
|D^2v(t,x) - D^2v(s,y)| \le m (1 + |x|^{p_5} + |y|^{p_5})(|t-s| + |x-y|)\\
\nonumber
\mbox{for all } 0 \le t,s \le T \mbox{ and } x,y \in
\mathbb{R}^d \, .
\eea
Note that \eqref{vlip} follows, for instance, if
$\frac{\partial}{\partial t} D^2 v$ and $D^3v$ exist and
$$
\left.
\begin{array}{l}
\max_{ij} |\frac{\partial}{\partial t} (D^2 v(t,x))_{ij} |\\
\max_{ij} |D (D^2 v(t,x))_{ij} |
\end{array}
\right\}
\le \frac{m}{d}(1 + |x|^{p_5}) \quad \mbox{for all } 0 \le t \le T \mbox{
and } x \in \mathbb{R}^d \, .
$$
Fix $(t,x) \in [0,T) \times \mathbb{R}^d$.
By Section \ref{sectPDE2BSDE}, the processes
\beas
Y_s &=& v(s,X^{t,x}_s) \, , \quad s \in [0,T] \, ,\\
Z_s &=& Dv(s, X^{t,x}_s) \, , \quad s \in [0,T] \, ,\\
\Gamma_s &=& D^2v(s, X^{t,x}_s) \, , \quad s \in [0,T] \, ,\\
A_s &=& {\cal L} Dv (s, X^{t,x}_s) \, , \quad s \in [0,T] \, ,
\eeas
solve the 2BSDE corresponding to $(X^{t,x}, f,g)$.
By \eqref{vbound} and \eqref{vlip}, $Z$ is in ${\cal A}^{t,x}_m$
(see \eqref{Zbound} and \eqref{gammalip}). Hence,
if the assumptions of Theorem \ref{thmmain} are fulfilled,
$(Y,Z, \Gamma, A)$ is the only solution of the 2BSDE corresponding to
$(X^{t,x}, f, g)$ with $Z \in {\cal A}^{t,x}$.
\end{Remark}

\noindent {\bf Proof of Theorem \ref{thmmain}} $\mbox{}$\\
{\bf 1.}
Let $(X^{0,x_0},Y^{0,x_0},\Gamma^{0,x_0},A^{0,x_0})$ be a solution
to the 2BSDE corresponding to $(X^{0,x_0},f,g)$
with $Z^{0,x_0} \in {\cal A}_m^{0,x_0}$ for some $m \ge 0$. Then, it follows from
Lemma \ref{pdpp} that
\be \label{YgeV}
Y^{0,x_0}_s \ge V^m(s,X^{0,x_0}_s) \quad \mbox{for all } s \in [0,T] \, ,
\ee
and by symmetry,
\be \label{YleU}
Y^{0,x_0}_s \le U_m(s,X^{0,x_0}_s) \quad
\mbox{for all } s \in [0,T] \, .
\ee
Recall that the inequalities \eqref{YgeV} and \eqref{YleU} are understood in the
$P$-almost sure sense. But since, by assumption,
$\sigma$ takes values in ${\cal M}_{\rm inv}$,
$X^{0,x_0}_s$ has full support for all $s \in (0,T]$
(see, for instance, Nualart, 1995), and we get
from \eqref{YgeV} and \eqref{YleU} that
$$
V^m(s,x) \le U_m(s,x) \quad \mbox{for all } (s,x) \mbox{ in a dense subset of }
[0,T] \times \mathbb{R}^d \, .
$$
It follows that
$$
V^m_* \le U_m^* \quad \mbox{on } [0,T] \times \mathbb{R}^d \, .
$$
Together with Proposition \ref{Vbound}, this shows that
$V^m_*$ and $U^*_m$ are finite on $[0,T] \times \mathbb{R}^d$.
By Remark \ref{viscsol},
$V^m_*$ is a viscosity supersolution and $U^*_m$ a viscosity subsolution of the PDE
\eqref{pde} on $[0,T) \times \mathbb{R}^d$.
Therefore, it follows from Proposition \ref{Vbound} and (A5) that
$$
V^m_* \ge U^*_m \quad \mbox{on } [0,T] \times \mathbb{R}^d \, .
$$
Hence, the function $v = V^m_* = V^m = U_m = U^*_m$ is continuous on
$[0,T] \times \mathbb{R}^d$ and
a viscosity solution with growth $p$ of the PDE \eqref{pde}--\eqref{terminal}.
By (A5), $v$ is the only viscosity solution of the PDE \eqref{pde}--\eqref{terminal}
with growth $p$. This shows 1.

\noindent {\bf 2}. Since $v = V^m = U_m$, it follows from
\eqref{YgeV} and \eqref{YleU} that
$$
Y^{0,x_0}_s = v(s, X^{0,x_0}_s) \quad \mbox{for all } s \in [0,T] \, .
$$
By conditioning $(Y^{0,x_0}, Z^{0,x_0}, \Gamma^{0,x_0}, A^{0,x_0})$
as in the proof of Lemma \ref{pdpp} and a
continuity argument, it can be deduced that for all $(t,x) \in
[0,T) \times \mathbb{R}^d$, there exists a solution
$(Y^{t,x}, Z^{t,x}, \Gamma^{t,x}, A^{t,x})$ of the 2BSDE corresponding to
$(X^{t,x},f,g)$ such that $Z^{t,x} \in {\cal A}^{t,x}$ and
$$
Y^{t,x}_s = v(s, X^{t,x}_s) \, , \quad s \in [t,T] \, .
$$
If for fixed $(t,x) \in [0,T) \times \mathbb{R}^d$,
$(Y,Z, \Gamma, A)$ is a solution of the 2BSDE corresponding to
$(X^{t,x},f,g)$ such that $Z \in {\cal A}^{t,x}$, then it follows as in 1. that
$$
Y_s = v(s, X^{t,x}_s) = Y^{t,x}_s \, , \quad s \in [t,T] \, ,
$$
and therefore also,
$$
(Z,\Gamma, A) = (Z^{t,x}, \Gamma^{t,x}, A^{t,x})
$$
because $Z^{t,x}$ is uniquely determined by
$$
Y^{t,x}_s = Y^{t,x}_t + \int_t^s f(r, X^{t,x_t}_r ,
Y^{t,x}_r, Z^{t,x}_r, \Gamma^{t,x}_r) dr
+ \int_t^s (Z^{t,x}_r)' \circ dX^{t,x}_r \, , \quad s \in [t,T] \, ,
$$
and $\Gamma^{t,x}$ and $A^{t,x}$ are uniquely determined by
$$
Z^{t,x}_s = Z^{t,x}_t + \int_t^s A^{t,x}_r dr +
\int_t^s \Gamma^{t,x}_r dX^{t,x}_r \, , \quad s \in [t,T] \, .
$$
This completes the proof of 2.
\unskip\nobreak\hfill$\Box$\par\addvspace{\medskipamount}

\setcounter{equation}{0}
\section{Monte Carlo methods for the solution of parabolic PDEs}
\label{sectnumerical}

In this section, we provide a formal discussion of the numerical
implications of our representation results. We start by recalling
some well-known facts in the linear case. We then review some
recent advances in the semi- and quasi-linear cases, and conclude with the fully
non-linear case related to Theorem \ref{thmmain}.

\subsection{The linear case}
\label{sectnumericallin}

In this subsection, we assume that the function $f$ is of the form
$$
f(t,x,y,z,\gamma) =
- \alpha(t,x) - \beta(t,x) y - \mu(x)' z
- \frac{1}{2} {\rm Tr} \left[\sigma(x) \sigma(x)' \gamma \right]
$$
Then, \eqref{pde} is a linear parabolic PDE. Under standard
conditions, it has a smooth solution $v$, and the Feynman-Kac representation theorem
states that for all $(t,x) \in [0,T] \times \mathbb{R}^d$,
$$
v(t,x) = \E{\int_t^T B_{t,s} \; \alpha \left(s,X^{t,x}_s \right) ds
+ B_{t,T} \; g \left(X^{t,x}_T \right)} \, ,
$$
where
$$
B_{t,s} := \exp \brak{\int_t^s \beta \left(r,X^{t,x}_r \right) dr}
$$
(see, for instance, Theorem 5.7.6 in Karatzas and Shreve, 1991).
This representation suggests a numerical approximation of the function $v$ by means of
the so-called Monte Carlo method:

(i) Given $J$ independent copies
$\left\{X^j ,~1\le j \le J \right\}$ of the process $X^{t,x}$, set
\beas
\hat v^{(J)}(t,x)
&:=&
\frac{1}{J} \sum_{j=1}^J \int_t^T B^j_{t,s}\; \alpha \left(s,X^j_s \right)ds
+ B^j_{t,T} \; g \left(X^j_T \right) \, ,
\eeas
where $B^j_{t,s} := \exp \brak{\int_t^s \beta \left(r,X^j_r \right)dr}$.
Then, it follows from the law of large numbers and the central limit theorem that
$$
\hat{v}^{(J)}(t,x) \to
v(t,x) \quad \mbox{a.s} \quad \mbox{and} \quad
\sqrt{J}\left(\hat v^{(J)}(t,x) - v(t,x) \right)
\to {\bf N}\left(0,\rho\right) \quad \mbox{in distribution} \, ,
$$
where $\rho$ is the variance of the random variable $\int_t^T
B_{t,s} \alpha \left(s,X^{t,x}_s \right)ds
+ B_{t,T}  \;g\left(X^{t,x}_T \right)$. Hence, $\hat{v}^{(J)}(t,x)$ is
a consistent approximation of $v(t,x)$.  Moreover, in contrast to
finite differences or finite elements methods, the error
estimate is of order $J^{-1/2}$, independently of the dimension $d$.

(ii) In practice, it is not possible to produce
independent copies $\left\{X^j,~1\le j \le J \right\}$
of the process $X^{t,x}$, except in trivial cases. In most
cases, the above Monte Carlo approximation is performed by
replacing the process $X^{t,x}$ by a suitable discrete-time
approximation $X^N$ with time step of order $N^{-1}$ for which
independent copies $\left\{X^{N,j},~1\le j \le J \right\}$ can be produced.
The simplest discrete-time approximation is the following
discrete Euler scheme: Set $X^N_t = x$ and for $1 \le n \le N$,
$$
X^N_{t_n} =  X^N_{t_{n-1}}
+ \mu(X^N_{t_{n-1}}) (t_n - t_{n-1})
+ \sigma (X^N_{t_{n-1}}) (W_{t_n} - W_{t_{n-1}}) \, ,
$$
where $t_n := t + n(T-t)/N$. We refer to Talay (1996) for a survey
of the main results in this area.

\subsection{The semi-linear case}
\label{sectnumericalsl}

We next consider the case where $f$ is given by
$$
f(t,x,y,z,\gamma) = \varphi(t,x,y,z)
- \mu(x)' z  - \frac{1}{2} \,{\rm  Tr}\left[\sigma(x) \sigma(x)'\gamma \right] \, .
$$
Then the PDE \eqref{pde} is called semi-linear. We assume that the
assumptions of Theorem \ref{thmmain} are satisfied.
In view of the connection between Fisk--Stratonovich and It\^{o}
integration, the 2BSDE \eqref{2bsde1}--\eqref{2bsde3}
reduces to an uncoupled FBSDE of the form
$$
\begin{aligned}
dY_s &= \varphi(s, X^{t,x}_s, Y_s , Z_s) ds
+ Z'_s \sigma (X^{t,x}_s) dW_s \, , \quad s \in [t,T) \, ,\\
Y_T &= g(X^{t,x}_T) \, ,
\end{aligned}
$$
(compare to Peng, 1991, 1992b; Pardoux and Peng 1992).
For $N \ge 1$, we denote
$t_n := t + n (T-t)/N$, $n = 0 , \dots, N$, and we define the discrete-time
approximation $Y^N$ of $Y$ by the backward scheme
$$
Y^N_T := g (X^{t,x}_T) \, ,
$$
and, for $n= 1, \ldots, N$,
\be\label{yn}
Y^N_{t_{n-1}} :=
\E{\left.Y^N_{t_n} \right| X^{t,x}_{t_{n-1}}}
- \varphi\left(t_{n-1}, X^{t,x}_{t_{n-1}},Y^N_{t_{n-1}},Z^N_{t_{n-1}}\right)
\left(t_n- t_{n-1}\right)
\ee
\be\label{zn}
Z^N_{t_{n-1}} :=
\frac{1}{t_n - t_{n-1}} \,
(\sigma (X^{t,x}_{t_{n-1}})')^{-1} \,
\E{(W_{t_n} - W_{t_{n-1}}) Y^N_{t_n} \mid X^{t,x}_{t_{n-1}}} \, .
\ee
Then, we have
$$
\limsup_{N \to \infty} \;
\sqrt{N} \;\left| Y^N_t - v(t,x) \right| < \infty \, ,
$$
see for instance, Zhang (2001), Bally and Pag\`es (2002),
Bouchard and Touzi (2004). The practical
implementation of this backward scheme requires the computation of
the conditional expectations appearing in \eqref{yn} and \eqref{zn}.
This suggests the use of a Monte Carlo approximation, as in the
linear case.  But at every time step, we need to compute
conditional expectations based on $J$ independent copies
$\left\{X^j,~1 \le j \le J \right\}$ of the process
$X^{t,x}$.  Recently, several approaches to this problem have
been developed. We refer to Lions and Regnier (2001), Bally and Pag\`es (2003),
Glasserman (2003), Bouchard and Touzi (2004) and the references
therein for the methodology and the analysis of such non-linear
Monte Carlo methods.

\subsection{The quasi-linear case}
\label{sectnumericalql}

It is shown in Antonelli (1993) and Ma et al. (1994) that coupled FBSDEs of
the form
$$
\begin{aligned}
dX_s &= \mu(s,X_s, Y_s, Z_s) ds + \sigma(s,X_s,Y_s) dW_s \, , \quad s \in [t,T] \, ,\\
X_t &= x\\
dY_s &= \varphi(s,X_s, Y_s, Z_s) ds + Z_s' \sigma(s,X_s,Y_s) dW_s \, ,
 \quad s \in [t,T) \, ,\\
Y_T &= g(X_T) \, ,
\end{aligned}
$$
are related to quasi-linear PDEs of the form \eqref{pde}--\eqref{terminal} with
$$
f(t,x,y,z,\gamma) = \varphi(t,x,y,z)
- \mu(t,x, y,z)' z  -
\frac{1}{2} \,{\rm  Tr}\left[\sigma(t,x,y) \sigma(t,x,y)'\gamma \right] \, ,
$$
see also, Pardoux and Tang (1999). Delarue and Menozzi (2004) have
used this relation to build a Monte Carlo scheme for the numerical
solution of quasi-linear parabolic PDEs.

\subsection{The fully non-linear case}
\label{sectnumericalnl}

We now discuss the case of a general $f$ as in the previous
section. Set
\bea
\varphi(t,x,y,z,\gamma) &:=& f(t,x,y,z,\gamma)
+ \mu(x)' z + \frac{1}{2} \,{\rm Tr} \left[\sigma(x) \sigma(x)' \gamma \right] \, .
\eea
Then for all $(t,x) \in [0,T) \times \mathbb{R}^d$ the 2BSDE
corresponding to $(X^{t,x}, f ,g)$ can be written as
\bea
\label{2bsdenumerical1}
dY_s &=& \varphi(s,X^{t,x}_s , Y_s, Z_s, \Gamma_s) ds + Z'_s
\sigma(X^{t,x}_s) dW_s \, , \quad s \in [t,T) \, ,\\
\label{2bsdenumerical2}
dZ_s &=& A_s ds + \Gamma_s dX^{t,x}_s \, , \quad s \in [t,T) \, ,\\
\label{2bsdenumerical3}
Y_T &=& g(X^{t,x}_T) \, .
\eea
We assume that the conditions of Theorem \ref{thmmain}
hold true, so that the PDE \eqref{pde} has a unique
viscosity solution $v$ with growth
$p = \max \crl{p_2 \, , \, p_3 \, , \, p_2 p_4 \, , \, p_4 + 2 p_1}$,
and there exists a unique
solution $(Y^{t,x},Z^{t,x},\Gamma^{t,x},A^{t,x})$ to the 2BSDE
\eqref{2bsdenumerical1}--\eqref{2bsdenumerical3}
with $Z^{t,x} \in {\cal A}^{t,x}$.

Comparing with the backward scheme
\eqref{yn}--\eqref{zn} in the semi-linear case, we suggest the
following discrete-time approximation of the processes
$Y^{t,x}$, $Z^{t,x}$ and $\Gamma^{t,x}$:
$$
Y^N_T := g (X^{t,x}_T) \, , \quad
Z^N_T := D g (X^{t,x}_T) \, ,
$$
and, for $n =1,\ldots,N$,
\beas
Y^N_{t_{n-1}} &:=& \E{\left.Y^N_{t_n} \right| X^{t,x}_{t_{n-1}}}
- \varphi (t_{n-1}, X^{t,x}_{t_{n-1}},Y^N_{t_{n-1}},Z^N_{t_{n-1}},
\Gamma^N_{t_{n-1}}) (t_n - t_{n-1})\\
Z^N_{t_{n-1}} &:=& \frac{1}{t_n - t_{n-1}} \,
(\sigma (X^{t,x}_{t_{n-1}})')^{-1} \,
\E{(W_{t_n} - W_{t_{n-1}}) Y^N_{t_n} \mid X^{t,x}_{t_{n-1}}} \, ,\\
\Gamma^N_{t_{n-1}} &:=& \frac{1}{t_n - t_{n-1}} \,
\E{Z^N_{t_n} \, (W_{t_n} - W_{t_{n-1}})' \mid X^{t,x}_{t_{n-1}}}
\sigma (X^{t,x}_{t_{n-1}})^{-1} \, .
\eeas
A precise analysis of the latter backward scheme is left for future
research. We conjecture that
$$
Y^N_t \to v(t,x) \quad \mbox{as} \quad N \to \infty \, ,
$$
if $v$ is $C^2$, then
$$
Z^N_t \to D v(t,x) \quad \mbox{as} \quad N \to \infty \, ,
$$
and if $v$ is $C^3$, then
$$
\Gamma^N_t \to D^2 v(t,x) \quad \mbox{as} \quad N \to \infty \, .
$$

\subsection{Connections with standard stochastic control}
\label{sectnumericalstoch}

Let ${\cal U}_0$ be the collection of all progressively measurable
processes $(\nu_t)_{t \in [0,T]}$ with values in a given bounded
subset $U \subset \mathbb{R}^k$.
Let $b : [0,T] \times \mathbb{R}^d \times U \to \mathbb{R}^d$ and
$a : [0,T] \times \mathbb{R}^d \times U \to {\cal S}^d$ be
continuous functions, Lipschitz in $x$ uniformly in
$(t,u)$. We call a process
$\nu \in {\cal U}_0$ an admissible control if
$$
\E{\int_0^T (|b(s,x,\nu_s)| + |a(s,x,\nu_s)|^2) ds} < \infty
$$
for all $x \in \mathbb{R}^d$ and denote the
class of all admissible controls by ${\cal U}$. For every pair of
initial conditions $(t,x) \in [0,T] \times \mathbb{R}^d$ and
each admissible control process $\nu \in {\cal U}$, the SDE
\be \label{Xnu}
\begin{aligned}
dX_s &= b(s,X_s,\nu_s) ds + a(s,X_s,\nu_s) dW_s \, , \quad s \in [t,T] \, ,\\
X_t &= x \, ,
\end{aligned}
\ee
has a unique strong solution, which we denote by
$(X^{t,x,\nu}_s)_{s \in [t,T]}$.
Let $\alpha, \beta : [0,T] \times \mathbb{R}^d \times U \to \mathbb{R}$ and
$g : \mathbb{R}^d \to \mathbb{R}$ be continuous functions with
$\beta \le 0$, and assume that $\alpha$ and $g$ have quadratic growth in
$x$ uniformly in $(t,u)$. Consider the stochastic control problem
$$
v(t,x) := \sup_{\nu \in {\cal U}}
\E{\int_t^T B^{\nu}_{t,s} \alpha(s,X^{t,x,\nu}_s, \nu_s) ds +
B^{\nu}_{t,T} \, g(X^{t,x,\nu}_T)} \, ,
$$
where
$$
B^{\nu}_{t,s} := \exp \brak{\int_t^s \beta(r,X^{t,x,\nu}_r,\nu_r) dr}
\, , \quad 0 \le t \le s \le T \, .
$$
By the classical method of dynamic programming, the function $v$ can be
shown to solve the Hamilton--Jacobi--Bellman equation
\be \label{HJB}
\begin{aligned}
&- v_t(t,x) + f(t,x,v(t,x), Dv(t,x), D^2v(t,x)) = 0\\
&  v(T,x) = g(x) \, ,
\end{aligned}
\ee
where
$$
f(t,x,y,z,\gamma) := \sup_{u \in U}
\crl{\alpha(t,x,u) + \beta(t,x,u) y + b(t,x,u)' z
- \frac{1}{2} {\rm Tr} [a a'(t,x,u) \, \gamma]} \, .
$$
This is a fully non-linear parabolic PDE covered by the class
\eqref{pde}--\eqref{terminal}. Note that $f$
is convex in the triple $(y,z,\gamma)$. The semi-linear case is
obtained when there is no control on the diffusion part, that is,
$a(t,x,u) = a(t,x)$ is independent of $u$.

If the value function $v$ has a stochastic representation in terms of a
2BSDE satisfying the assumptions of Theorem \eqref{thmmain},
then the Monte Carlo scheme of the
previous subsection can be applied to approximate $v$.

Under suitable regularity assumptions, the optimal control at time $t$
is known to be of the form
$$
\hat{u}(t,x,v(t,x), Dv(t,x), D^2 v(t,x)) \, ,
$$
where $\hat{u}$ is a maximizer of the expression
$$
\sup_{u \in U} \crl{\alpha(t,x,u) + \beta(t,x,u) v(t,x) + b(t,x,u)' Dv(t,x)
- \frac{1}{2} {\rm Tr} [a a'(t,x,u) D^2 v(t,x)]} \, .
$$

Notice that the numerical scheme suggested in
the Subsection \ref{sectnumericalnl} calculates at each step in time the values of
the processes $(X^{t,x},Y^{t,x},Z^{t,x},\Gamma^{t,x})$. Therefore, the optimal control is
also provided by this numerical scheme by
$$
\hat{\nu}_s = \hat{u} (s, X^{t,x}_s , Y^{t,x}_s, Z^{t,x}_s, \Gamma^{t,x}_s) \, .
$$

\setcounter{equation}{0}
\section{Boundary value problems}

In this section, we give a brief outline of an extension to
boundary value problems. Namely, let $O \subset \mathbb{R}^d$ be an open
set. For $(t,x) \in [0,T) \times \mathbb{R}^d$, the process
$X^{t,x}$ is given as in \eqref{sde}, but we stop it at the boundary
of $O$. Then we extend the terminal condition \eqref{2bsde3} in the 2BSDE to
a boundary condition. In other words, we introduce the exit time
$$
\theta := \inf \crl{s \ge t \mid X^{t,x}_s \notin O}
$$
and modify the 2BSDE \eqref{2bsde1}--\eqref{2bsde3} to
$$
\begin{aligned}
& Y_{s \wedge \theta} = g(X^{t,x}_{T \wedge \theta})
- \int_{s \wedge \theta}^{T \wedge \theta}
f(r,X^{t,x}_r,Y_r, \Gamma_r, A_r) dr -
\int_{s \wedge \theta}^{T \wedge \theta} Z'_r \circ dX^{t,x}_r \, ,
\quad s \in [t,T] \, ,\\
& Z_{s \wedge \theta} = Z_{T \wedge \theta} -
\int_{s \wedge \theta}^{T \wedge \theta} A_r dr
- \int_{s \wedge \theta}^{T \wedge \theta} \Gamma_r dX^{t,x}_r \, ,
\quad s \in [t,T] \, .
\end{aligned}
$$
Then the corresponding PDE is the same as \eqref{pde}--\eqref{terminal},
but it only holds
in $[0,T) \times O$. Also, the terminal condition $v(T,x) = g(x)$
only holds in $O$. In addition, the following lateral boundary condition holds
$$
v(t,x) = g(x) \, , \quad \mbox{for all } (t,x) \in [0,T] \times \partial O \, .
$$
All the results of the previous sections can easily be adapted to this case. Moreover,
if we assume that $O$ is bounded, most of the technicalities related to
the growth of solutions are avoided as the solutions are expected
to be bounded.


\bigskip \bigskip
\noindent
\begin{thebibliography}{aa12}

\bibitem{ant}
Antonelli, F. (1993).  Backward-forward stochastic differential
equations. {\sl Annals of Applied Probability} 3, 777--793.

\bibitem{as}
Arkin, V. and Saksonov, M. (1979). Necessary optimality conditions
for stochastic differential equations. {\sl Soviet Math. Dokl.}
20, 1--5.

\bibitem{bapa}
Bally, V., Pag\`es, G. (2003). Error analysis of the quantization
algorithm for obstacle problems. {\sl Stochastic Processes and
their Applications} 106(1), 1--40.

\bibitem{bbbl}
Barles, G., Biton, S., Bourgoing, M., Ley, O. (2003).
Uniqueness results for quasilinear parabolic equations through
viscosity solutions' methods.  {\sl Calc. Var. PDEs}, 18
(2),159-179.

\bibitem{bbp}
Barles, G., Buckdahn, R., and Pardoux, E. (1997). Backward
stochastic differential equations and integral-partial
differential equations, {\sl Stochastics Stochastics Rep.}, 60,
no. 1-2, 57--83.

\bibitem{bp}
Barles, G. and Perthame, B. (1988). Exit time problems in optimal
control and vanishing viscosity solutions of Hamilton--Jacobi
equations, {\sl SIAM J. Cont. Opt.} 26, 1133--1148.

\bibitem{bi73}
Bismut, J.M. (1973). Conjugate convex functions in optimal
stochastic control, {\sl J. Math. Anal. Appl.} 44, 384--404.

\bibitem{bi78}
Bismut, J.M. (1978). Contr\^{o}le des syst\`eme lin\'eaire
quadratiques: applications de l'integrale stochastique, {\sl
S\'em. Probab.} XII, Lect. Notes. Math.
vol. 649. Springer-Velag, 180--264.

\bibitem{bt}
Bouchard, B., Touzi, N. (2004). Discrete-time approximation and
Monte Carlo simulation of backward stochastic differential
equations, {\sl Stochastic Processes and their Applications} 111,
175--206.

\bibitem{caca}
Cabre, X., Caffarelli, L. (1995). {\sl Fully Nonlinear
Elliptic Equations}, AMS, Providence.

\bibitem{csta}
Cheridito, P., Soner, H.M., Touzi, N. (2005a). Small time path
behavior of double stochastic integrals and applications to
stochastic control, {\sl Annals of Applied Probability}, to appear.

\bibitem{cstb}
Cheridito, P., Soner, H.M., Touzi, N. (2005b). The
multi-dimensional super-replication problem under gamma
constraints, {\sl Annales de l'Institute Henri Poincar\'e (C)
Non Linear Analysis}, to appear.

\bibitem{chevance}
Chevance, D. (1997). Numerical Methods for Backward Stochastic
Differential Equations, {\sl Publ. Newton Inst.}, Cambridge University Press.

\bibitem{cil}
Crandall, M.G., Ishii, H., Lions, P.L. (1992). User's guide to
viscosity solutions of second order partial differential
equations. {\sl Bull. Amer. Math. Soc.} 27(1), 1--67.

\bibitem{cma}
Cvitani\'c, J., Ma, J. (1996). Hedging options for a large
investor and forward-backward SDE's, {\sl Ann. Appl. Probab.} 6,
370--398.

\bibitem{cvk}
Cvitani\'c, J., Karatzas, I. (1996).  Backward stochastic
differential equations with reflection and Dynkin games. {\sl
Annals of Probability}, 24(4), 2024--2056.

\bibitem{cks}
Cvitani\'c, J., Karatzas, I., Soner, H.M. (1999). Backward SDEs
with constraints on the gains-process, {\sl Annals of
Probability} 26, 1522--1551.

\bibitem{delarue}
Delarue, F. (2002). On the existence and uniqueness of solutions
to FBSDEs in a non-degenerate case, {\sl Stoch. Proc. Appl.} 99,
209--286.

\bibitem{demaz}
Delarue, F., Menozzi, F. (2004). A forward-backward stochastic
algorithm for quasi-linear PDEs. Preprint 932, University of Paris VI \& VII.

\bibitem{dmap}
Douglas, J., Ma, J., Protter, P. (1996). Numerical methods for
forward-backward stochastic differential equations. {\sl Ann. Appl. Probab.}
6, 940--968.

\bibitem{elpq}
El Karoui, N., Peng, S., Quenez, M.C. (1997).  Backward
stochastic differential equations in finance.  {\sl Mathematical
Finance} 7(1), 1--71.

\bibitem{evans}
Evans, L.C. (1998). {\sl Partial Differential Equations}. AMS,
Providence.

\bibitem{fs}
Fleming, W.H., Soner, H.M. (1993). {\sl Controlled Markov
Processes and Viscosity Solutions}. Applications of Mathematics
25. Springer-Verlag, New York.

\bibitem{Glasserman}
Glasserman, P. (2003). {\sl Monte Carlo Methods in Financial
Engineering}, Stochastic Modelling and Applied Probability, Vol.
53. Springer.

\bibitem{ikwa}
Ikeda, N., Watanabe, S. (1989). {\sl Stochastic Differential
Equations and Diffusion Processes}, Second Edition.
North-Holland Publishing Company.

\bibitem{ishii}
Ishii, H. (1984). Uniqueness of unbounded viscosity solutions of
Hamilton-Jacobi equations. {\sl Indiana U. Math. J.}. 33,
721--748.

\bibitem{ks91}
Karatzas, I., Shreve, S. (1991). {\sl Brownian Motion and Stochastic Calculus},
Second Edition. Springer-Verlag.

\bibitem{krylov}
Krylov, N. (1987). {\sl Nonlinear elliptic and Parabolic Partial
Differential Equations of Second Order}. Mathematics and its
Applications, Reider.

\bibitem{lsu}
Ladyzenskaya, O.A., Solonnikov, V.A., Uralseva, N.N. (1967). {\sl
Linear and Quasilinear Equations of Parabolic Type}, AMS,
Providence.

\bibitem{lionsr}
Lions, P.L., Regnier, H. (2001). Calcul du prix et des
sensibilit\'es d'une option am\'ericaine par une m\'ethode de
Monte Carlo. Preprint.

\bibitem{map}
Ma J., Protter, P., Yong, J. (1994). Solving backward stochastic
differential equations explicitley - A four step
scheme. {\sl Prob. Theory Rel. Fields} 98, 339--359.

\bibitem{mayong}
Ma, J., Yong, J. (1999). {\sl Forward-Backward Stochastic
Differential Equations and their Applications}.
LNM 1702, Springer-Verlag.

\bibitem{mapsmt}
Ma, J., Protter, P., San Martin, J., Torres, S. (2002). Numerical
methods for  backward stochastic differential equations,
{\sl Ann. Appl. Probab.} 12, 302--316.

\bibitem{nualart}
Nualart, D. (1995). {\sl The Malliavin Calculus and Related
Topics}. Springer-Verlag.

\bibitem{pp90}
Pardoux, E., Peng, S. (1990). Adapted solution of a backward
stochastic differential equation, {\sl Systems Control Lett.}, 14, 55--61.

\bibitem{pp92}
Pardoux, E., Peng, S. (1992a). {\sl Backward stochastic
differential equations and quasilinear parabolic partial
differential equations.} Lecture Notes in CIS, Vol. 176.
Springer-Verlag, 200--217.

\bibitem{pp94}
Pardoux, E., Peng, S. (1994). Backward doubly stochastic
differential equations and systems of quasilinear parabolic
SPDEs, {\sl Probab. Theory Rel. Fields} 98, 209--227.

\bibitem{pardouxtang}
Pardoux, E., Tang, S. (1999). Forward-backward stochastic
differential equations and quasilinear parabolic
PDEs. {\sl Prob. Th. Rel. Fields} 114(2), 123--150.

\bibitem{peng90}
Peng, S. (1990). A general stochastic maximum principle for
optimal control problems, {\sl SIAM J. Control Optim.} 28,
966--979.

\bibitem{peng91}
Peng, S. (1991). Probabilistic interpretation for systems of
quasilinear parabolic partial differential equations, {\sl Stochastics}
37, 61--74.

\bibitem{peng92a}
Peng, S. (1992a) Stochastic Hamilton--Jacobi--Bellman equations,
{\sl SIAM J. Control Optim.} 30, 284--304.

\bibitem{peng92b}
Peng, S. (1992b). A generalized dynamic programming principle and
Hamilton--Jacobi--Bellman equation, {\sl Stochastics} 38,
119--134.

\bibitem{peng92c}
Peng, S. (1992c). A nonlinear Feynman--Kac formula and
applications, in {\sl Proccedings of Symposium of System Sciences
and Control Theory,}
Ed. Chen and Yong. Singapore: World Scientific, 173--184.

\bibitem{peng93}
Peng, S. (1993). Backward stochastic differential equation and
its application in optimal control, {\sl Appl. Math. Optim.} 27,
125--144.

\bibitem{sonertouzi}
Soner, H.M., Touzi, N. (2002). Stochastic target
problems, dynamic programming and viscosity solutions, {\sl SIAM
Journal on Control and Optimization} 41, 404--424.

\bibitem{talay}
Talay, D. (1996). Probabilistic numerical methods for partial
differential equations: elements of analysis, in {\sl
Probabilistic Models for Nonlinear Partial Differential
Equations}, D. Talay and L. Tubaro, editors, Lecture Notes in
Mathematics 1627,
48--196, 1996.

\bibitem{zhang}
Zhang, J. (2001). Some fine properties of backward stochastic
differential equations. PhD thesis, Purdue University.

\end{thebibliography}
\end{document}